\documentstyle[12pt,fleqn]{article}

\font\twlgot =eufm10 scaled \magstep1
\font\egtgot =eufm8
\font\sevgot =eufm7
\font\twlmsb =msbm10 scaled \magstep1
\font\egtmsb =msbm8
\font\sevmsb =msbm7

\newfam\gotfam

\textfont\gotfam\twlgot
\scriptfont\gotfam\egtgot
\scriptscriptfont\gotfam\sevgot

\newfam\msbfam
\textfont\msbfam\twlmsb
\scriptfont\msbfam\egtmsb
\scriptscriptfont\msbfam\sevmsb
\def\Bbb{\protect\pBbb}
\def\pBbb{\relax\ifmmode\expandafter\Bb\else\typeout{You cann't use
Bbb in text mode}\fi}
\def\Bb #1{{\fam\msbfam\relax#1}}

\def\thebibliography#1{\section*{References}\list
  {[\arabic{enumi}]}{\settowidth\labelwidth{#1}\leftmargin\labelwidth
    \advance\leftmargin\labelsep
    \usecounter{enumi}}
    \def\newblock{\hskip .11em plus .33em minus .07em}
    \sloppy\clubpenalty4000\widowpenalty4000
    \sfcode`\.=1000\relax}

\def\op#1{\mathop{\fam0 #1}\limits}

\newcommand{\beq}{\begin{equation}}
\newcommand{\eeq}{\end{equation}}
\newcommand{\ben}{\begin{eqnarray}}
\newcommand{\een}{\end{eqnarray}}
\newcommand{\be}{\begin{eqnarray*}}
\newcommand{\ee}{\end{eqnarray*}}
\newcommand{\bea}{\begin{eqalph}}
\newcommand{\eea}{\end{eqalph}}

\newcommand{\cH}{{\cal H}}

\newcommand{\la}{\lambda}

\newcommand{\f}{\phi}
\newcommand{\om}{\omega}
\newcommand{\Om}{\Omega}
\newcommand{\m}{\mu}
\newcommand{\n}{\nu}

\newcommand{\vt}{\vartheta}
\newcommand{\vf}{\varphi}

\newcommand{\w}{\wedge}

\newcommand{\dr}{\partial}

\newcommand{\fl}{\flat}
\newcommand{\sh}{\sharp}

\newcounter{eqalph}
\newcounter{equationa}
\newcounter{theorem}
\newcounter{remark}
\newcounter{proposition}
\newcounter{lemma}
\newcounter{corollary}
\newcounter{definition}

\newenvironment{eqalph}{\stepcounter{equation}
\setcounter{equationa}{\value{equation}}
\setcounter{equation}{0}

\begin{eqnarray}}{\end{eqnarray}\setcounter{equation}{\value{equationa}}}

\def\theremark{\arabic{remark}}

\def\thedefinition{\arabic{definition}}

\newcommand{\mar}[1]{}

\begin{document}
\hbox{}

{\parindent=0pt

{\large\bf A note on the KAM theorem for partially integrable Hamiltonian systems}
\bigskip

{\sc G. GIACHETTA}$^1$, {\sc L. MANGIAROTTI}$^2$ and {\sc G.SARDANASHVILY}$^3$
\bigskip

\begin{small}

$^1${\it Department of Mathematics and Informatics, University of Camerino, 62032
Camerino (MC), Italy 

E-mail: Giovanni.Giachetta@unicam.it}

$^2${\it Department of Mathematics and Informatics, University of Camerino, 62032
Camerino (MC), Italy

E-mail: Luigi.Mangiarotti@unicam.it} 

$^3${\it Department of Theoretical Physics, Physics Faculty, Moscow State
University, 117234 Moscow, Russia

 E-mail sard@grav.phys.msu.su, 
sard@campus.unicam.it, URL: http://webcenter.ru/$\sim$sardan/}
\bigskip

{\bf Abstract.} We provide a symplectic reduction of a partially integrable 
Hamiltonian system
to a completely integrable one. The KAM theorem is applied to
this reduced completely integrable Hamiltonian system. Its KAM perturbation
generates a perturbation of the original partially integrable system
such that 
the set of tori kept by this perturbation is large. 

\medskip

\medskip


\end{small}
}

\bigskip
\bigskip

Without loss of generality, we here
deal with the KAM theorem in its standard form of quasi-periodic stability and 
anyliticity \cite{arn,ben,laz}.
Let
\mar{k61}\beq
V\times T^m\to V \label{k61}
\eeq
be a trivial fibre bundle in $m$-dimensional tori over an $m$-dimensional domain $V\subset \Bbb R^m$.
Let it be equipped with the standard coordinates $(J_a,\vf^a)$, $a=1,\ldots,m$,
and provided with the symplectic form
\mar{dd26}\beq
\Om= dJ_a\w d\vf^a. \label{dd26} 
\eeq
Let a Hamiltonian $\cH(J_a)$ on $V\times T^m$ depend
only on action variables $(J_a)$. Its Hamiltonian vector field with respect to the
symplectic form $\Om$ (\ref{dd26}) reads
\mar{bi49}\beq
\vt_\cH=\dr^a\cH\dr_a, \label{bi49}
\eeq
and defines the Hamilton equation
\mar{k60}\beq
\dot J_a=0, \qquad \dot\vf^a=\dr^a\cH \label{k60}
\eeq
on $V\times T^m$. By this equation, the tori of the fibre bundle (\ref{k61}) are invariant manifolds.
In accordance with the classical  Liouville--Arnold theorem \cite{arn,laz}, any
completely integrable Hamiltonian system around a regular connected compact invariant manifold
is brought into the above mentioned standard form. 
Now, let 
\mar{k12}\beq
\cH'=\cH+\cH_1(J_a,\vf^a) \label{k12}
\eeq
be a perturbed Hamiltonian on $V\times T^m$
depending both on action and angle coordinates. Its Hamiltonian vector
field with respect to the symplectic form $\Om$ (\ref{dd26}) reads
\mar{bi51}\beq
\vt_{\cH'}=\dr^a\cH'\dr_a -\dr_a\cH'\dr^a, \label{bi51}
\eeq
and defines the Hamilton equation
\mar{k62}\beq
\dot J_a=-\dr_a\cH', \qquad \dot\vf^a=\dr^a\cH' \label{k62}
\eeq
on $V\times T^m$. The above mentioned KAM theorem states 
that, under certain conditions on the Hamiltonian $\cH'$ (\ref{k12}), 
the perturbed Hamilton equation (\ref{k62})
admits solution living in tori so that the complement of
the set of these tori has a Lebesgue measure which tends
to zero as the perturbation tends to zero. 
 Namely, fixing a torus $J^*=(J^*_a)$
with appropriate (Diophantine)
frequencies $\om^a=(\dr^a\cH)(J^*_b)$, one can construct
a convergent sequence
of canonical analytical transformations $(J_a,\vf^a)=\Phi(p_a,q^a)$
around $J^*$ which brings the perturbed Hamiltonian $\cH'$ (\ref{k12})
into the form
\be
\cH'(p_a,q^a)= {\rm const.}+ \om^ap_a + \cH_2(p_a,q^a), 
\ee
where $\cH_2$ is a strictly quadratic function of canonical coordinates $p_a$. 
Then, the Hamilton equation 
(\ref{k62}),
rewritten in coordinates $(p_a,q^a)$, has a particular
solution 
\mar{k63}\beq
p_a=0, \qquad  q^a=\om^at \label{k63}
\eeq
living in a torus. In other words, the Hamiltonian vector field (\ref{bi51})
has the integral curve $\Phi(0,\om^at)$ (\ref{k63}) located in a torus. 
In terms of the original
action-angle coordinates, this torus is given by the parametric equation 
$(J_a,\vf^a)=\Phi(0,q^a)$, $q\in T^m$. It is an image of the 
torus $J^*$ under a diffeomorphism near the identity.

Turn now to a partially integrable Hamiltonian system 
on a $2n$-dimensional 
symplectic manifold $(Z,\Om)$. It is defined by $1\leq k\leq n$ 
 real smooth functions $\{H_i\}$ which are pairwise in involution 
and independent almost everywhere on $Z$. The latter implies that 
the set of non-regular points, where
the morphism
\mar{d1}\beq
\op\times^k H_i : Z\to\Bbb R^k \label{d1}
\eeq
fails to be a submersion, is nowhere dense. 
We agree to think of one of the functions $H_i$ as being a Hamiltonian
of a partially integrable system. 
Let $M$ be a regular connected
compact invariant manifold of this system.
The well-known Nekhoroshev theorem 
\cite{arn,gaeta,nekh} states that, under certain conditions, 
there exists an 
open neighbourhood $U$ of $M$ which is a trivial composite bundle 
\mar{z10}\beq
\pi:U=V\times W\times T^k\to V\times W\to V \label{z10}
\eeq
over domains $W\subset \Bbb R^{2(n-k)}$ and $V\subset \Bbb R^k$. It
is provided with the partial
action-angle coordinates 
$(I_i;z^\la; \f^i)$, $i=1,\ldots,k$, $\la=1,\ldots,2(n-k)$, such that 
the symplectic form $\Om$ on $U$ reads
\mar{d26}\beq
\Om= dI_i\w d\f^i +\Om_{\m\n}(I_j,z^\la) dz^\m\w dz^\n +\Om_\m^i(I_j,z^\la) dI_i\w dz^\m,
\label{d26} 
\eeq
and integrals of motion $H_i$ depend only on the action 
coordinates $I_j$. 
Let $\cH(I_i)$ be a Hamiltonian of a partially integrable system on $U$ (\ref{z10}). 
Its Hamiltonian vector field 
with respect to the symplectic form $\Om$ (\ref{d26}) is
\mar{k7}\beq
\vt_\cH=\dr^i\cH\dr_i, \label{k7}
\eeq
and defines the Hamilton equation
\mar{k65}\beq
\dot I_i=0, \qquad \dot z^\la=0, \qquad \dot\f^i=\dr^i\cH. \label{k65}
\eeq
A glance at this equation shows that $U\to V\times W$ (\ref{z10})
is a trivial bundle in invariant tori. Herewith, the functions $H^\la=z^\la$ on $U$ (\ref{z10})
are also integrals of motion, but they need not be in involution with the functions
$H_i$ because of the third term $\Om_\m^i dI_i\w dz^\m$ of the symplectic form
$\Om$ (\ref{d26}). If this term vanishes, there are more than $k$ integrals of
motion on $U$ in involution. 

Let us note that there exists a Darboux coordinate
chart $Q\times T^k\subset U$, foliated in tori and provided with 
coordinates $(I_i;p_\la;q^\la;\vf^i)$ where $\vf^i=\f^i + f^i(I_i, z^\la)$ are 
new coordinates on tori \cite{fasso}. The symplectic form $\Om$ (\ref{d26}) 
on this chart takes the canonical form
\be
\Om= dI_i\w d\vf^i + dp_\la\w dq^\la.
\ee
It follows that a partially integrable system in question on this 
chart becomes completely
integrable, but its invariant manifolds fail to be compact and
the KAM theorem is not applied to this completely integrable
system.

Now, let us consider a perturbed Hamiltonian 
\mar{k12'}\beq
\cH'=\cH+\cH_1(I_i,\f^i) \label{k12'}
\eeq
on  $U$ (\ref{z10}) depending both on action and angle variables.
Its Hamiltonian vector field with respect to the symplectic form $\Om$
(\ref{d26}) reads
\mar{k66}\beq
\vt_{\cH'}= (C^{ik}\dr_k\cH' +\dr^i\cH')\dr_i -\dr_i\cH'\dr^i
+ B^{\la k}\dr_k\cH'\dr_\la, \label{k66}
\eeq
where the coefficients $C^{ik}$ and $B^{\la k}$ are expressed into the components 
$\Om_{\m\n}$ and $\Om_\m^i$ of the symplectic form (\ref{d26}), and they vanish if
$\Om_\m^i$ do so. A comparison of the expressions (\ref{bi51}) and (\ref{k66})
shows that the Hamilton equation defined by the 
vector field (\ref{k66}) is not that we want. One can try to study non-Hamiltonian
perturbations of the Hamilton equation (\ref{k65}). However, this is not 
the case of KAM theorem because, as was mentioned above, the proof of this theorem
is based on the use of canonical transformations in order to show that
a perturbed equation has a solution living in a tori. 
Therefore, we will replace the symplectic structure (\ref{d26})
on $U$ with a different Poisson one.

Let us provide $U$ (\ref{z10}) with the degenerate Poisson bivector field
\mar{k0}\beq
w=\dr^i\w \dr_i \label{k0}
\eeq
of constant rank $k$. The corresponding Poisson bracket reads
\mar{k2}\beq
\{f,f'\}=\dr^if\dr_i f' - \dr_if\dr^if', \qquad f,f'\in C^\infty(U). 
\label{k2}
\eeq
It is readily observed that 
the functions $H_i$ are in involution with respect to this Poisson bracket. Moreover,
their Hamiltonian vector fields with respect both to the symplectic form
$\Om$ (\ref{d26}) and the Poisson bivector field $w$ (\ref{k0})
coincide. In particular, the Hamiltonian vector field of a non-perturbed
Hamiltonian $\cH$ relative to $w$ is also $\vt_\cH$ (\ref{k7}), and it defines
the same Hamilton equation (\ref{k65}) on $U$. Consequently, functions $H_i$
remain integrals of motion. It follows that the symplectic form
$\Om$ (\ref{d26}) and the Poisson bivector field $w$ (\ref{k0}) make up
a bi-Hamiltonian system. The corresponding recursion operator is written as
$R=w^\sh\circ \Om^\fl$ in terms of the morphisms $w^\sh:T^*U\to TU$ and 
$\Om^\fl: TU\to T^*U$.  

Furthermore, 
let us consider the alternative fibration
\mar{k1}\beq
\pi':U=V\times W\times T^k\to V\times T^k \label{k1}
\eeq
of the trivial bundle (\ref{z10}). Then, 
the Poisson bracket (\ref{k2})
on $U$ yields a Poisson bracket $\{,\}'$ on the toroidal cylinder 
$V\times T^k$ by the rule
\mar{k6}\beq
\pi'^*\{f,f'\}'=\{\pi'^*f,\pi'^*f'\}, 
\qquad f,f'\in C^\infty(V\times T^k). \label{k6}
\eeq
With respect to the action-angle coordinates $(I_i,\f^i)$ on $V\times T^k$,
the Poisson bracket $\{,\}'$ and the corresponding Poisson bivector field $w'$
on $V\times T^k$ 
have the same coordinate expressions (\ref{k2}) and (\ref{k1}), respectively,
as the bracket $\{,\}$  and the Poisson bivector field $w$ 
 on $U$. Herewith,
the Poisson bivector field $w'$ is non-degenerate and associated to the symplectic form
\mar{k4}\beq
\Om'=dI_i\w d\f^i \label{k4}
\eeq
on $V\times T^k$. Let $s$ be a section of the fibre bundle $\pi'$ (\ref{k1}) such that
$z^\la\circ s=$const. Then, one can think of the symplectic
manifold $(V\times T^k,\Om')$ as being
a reduction of the symplectic manifold $(U,\Om)$ via the submanifold
$s(V\times T^k)$. Furthermore, the pull-back functions $s^*H_i$,
denoted $H_i$ again, make up a completely integrable system on 
the symplectic manifold $(V\times T^k,\Om')$. The following property
enables us to treat this system as 
a reduction of the partially integrable system $\{H_i\}$ on $U$.
Due to the relation (\ref{k6}),
 any Hamiltonian vector field 
$\vt'=-w'\lfloor df$ on the symplectic manifold
$(V\times T^k,\Om')$ gives rise to the Hamiltonian vector field
\mar{k10}\beq
\vt=-w\lfloor (\pi'^*df) \label{k10}
\eeq
 on the Poisson manifold $(U,w)$. 
Moreover, the Hamiltonian vector fields of integrals of motion
$H_i$ on $V\times T^k$ give rise exactly to the Hamiltonian vector fields
of integrals of motion $H_i$ on the 
Poisson manifold $(U,w)$ which, as was mentioned above,
 coincide with Hamiltonian vector fields of
$H_i$ on the symplectic manifold $(U,\Om)$.

In particular, let $\cH(I_i)$ be the pull-back onto $V\times T^k$ 
of a non-perturbed Hamiltonian of a partially integrable system
on $U$ (\ref{z10}). One can think of it as being a Hamiltonian
of the reduced completely integrable system on $V\times T^k$.
Its Hamiltonian vector field takes the coordinate form 
(\ref{k7}) and yields the Hamilton equation
\mar{k67}\beq
\dot I_i=0, \qquad \dot \f^i=\dr^i\cH \label{k67}
\eeq
on $V\times T^k$ (cf. (\ref{k60})). Let $\cH'$ be a perturbation of the
Hamiltonian $\cH$ which fulfils the conditions of the KAM 
theorem on the symplectic manifold $(V\times T^k,\Om')$.
Then, its Hamiltonian vector field 
\mar{bi51'}\beq
\vt'_{\cH'}=\dr^i\cH'\dr_i -\dr_i\cH'\dr^i \label{bi51'}
\eeq
(cf. (\ref{bi51})) on $V\times T^k$ has an integral curve 
\be
I_i=\Phi_i(0,w^jt), \qquad \f^i=\Phi^i(0,w^jt)
\ee
located in a torus. By the formula (\ref{k10}), the Hamiltonian vector
field (\ref{bi51'}) gives rise
to the vector field
\mar{k68}\beq
\vt_{\cH'}=-w\lfloor (\pi'^*\cH')=\dr^i\cH'\dr_i -\dr_i\cH'\dr^i
\label{k68}
\eeq
on $U$. It defines the first order dynamic equation 
\mar{k11}\beq
\dot I_i=-\dr_i\cH', \qquad \dot z^\la=0, 
\qquad \dot \f^i=\dr^i\cH' \label{k11}
\eeq
on $U$. 
The vector field (\ref{k68}) has integral curves
\be
I_i=\Phi_i(0,w^jt), \qquad z^\la={\rm const.},
\qquad \f^i=\Phi^i(0,w^jt)
\ee
located in tori. These curves provide particular solutions of the
dynamic equation (\ref{k11}).
As follows from the expression (\ref{k68}), 
this equation is a Hamilton equation of the Hamiltonian $\pi'^*\cH'$
with respect to the 
Poisson structure $w$ (\ref{k0}) on $U$, but need not
be so relative to the original symplectic form (\ref{d26}).
One can think of (\ref{k11}) as being {\it sui generis} a 
KAM perturbation of 
the Hamilton equation (\ref{k65}). Clearly,
the set of tori kept under this perturbation on a fibre
$z^\la=$const. is as large as that guaranteed by the
KAM theorem on $V\times T^k$.

Basing on this result, one can extend the KAM theorem 
on a symplectic manifold $V\times T^m$ (\ref{k61}), $\Om$ 
(\ref{dd26}) 
to a Poisson manifold $U$ (\ref{z10}), $w$ (\ref{k0}) 
in a straightforward manner. It's proof is a repetition of that
on a symplectic manifold (e.g., in \cite{ben}), but appeals to 
canonical transformations of the degenerate Poisson structure
$w$ (\ref{k10}). 

Let us point out possible applications of this extended KAM theorem. 
These are completely integrable systems whose invariant manifold
need not be compact, but is foliated in tori. For instance, this is 
the case of time-dependent completely integrable systems \cite{acang1,acang2}.
Note that any closed mechanical system (e.g., an $n$-body system) in an 
Euclidean three-space has at least 3 integrals of motion, besides
a Hamiltonian, which are the components $P_1$, $P_2$, $P_3$
of its total momentum. If this system admits other integrals of
motion $H_i$ 
(which need not  be in involution with the total momentum $P_\la$), 
it is also the case of application
of the extended KAM theorem.

\end{document}